\documentclass[review]{elsarticle}

\journal{Journal of \LaTeX\ Templates}

\usepackage{amssymb,latexsym,amsmath,amsthm,graphicx}
\usepackage{geometry,amsfonts,amsthm}
\usepackage{color}
\usepackage[english]{babel}
\usepackage[T1]{fontenc}

\newcommand{\R}{{\mathbb R}}

\newcommand{\cald}{{\mathcal D}}
\newcommand{\calf}{{\mathcal F}}

\newcommand{\calu}{{\mathcal U}}

\newcommand{\cl}{\operatorname{cl}}

\newcommand{\dom}{\operatorname{dom}}

\newcommand{\fix}{\operatorname{fix}}
\newcommand{\gph}{\operatorname{gph}}
\newcommand{\inte}{\operatorname{int}}
\newcommand{\dist}{\operatorname{dist}}
\newcommand{\argmin}{\operatornamewithlimits{arg\,min}}

\newtheorem{Theorem}{Theorem}

\newtheorem{Corollary}{Corollary}[section]
\newtheorem{DEF}{Definition}[section]
\newenvironment{Definition}{\begin{DEF}\rm}{\end{DEF}}
\newtheorem{EX}{Example}[section]

\newtheorem{REM}{Remark}[section]
\newenvironment{Remark}{\begin{REM}\rm}{\end{REM}}

\begin{document}

\begin{frontmatter}

\title{A continuity result for the adjusted normal cone operator}



\author[1]{Marco Castellani}
\ead{marco.castellani@univaq.it}

\author[1]{Massimiliano Giuli}
\ead{massimiliano.giuli@univaq.it}


\address[1]{Department of Information Engineering, Computer Science and Mathematics, University of L'Aquila, Via Vetoio, L'Aquila, Italy}

\begin{abstract}
The concept of adjusted sublevel set for a quasiconvex function was introduced in \cite{AuHa05} and the local existence of a norm-to-weak$^*$ upper semicontinuous base-valued submap of the normal operator associated to the adjusted sublevel set was proved.
When the space is finite dimensional, a globally defined upper semicontinuous base-valued submap is obtained taking the intersection of the unit sphere, which is compact, with the normal operator, which is closed.
Unfortunately, this technique does not work in the infinite dimensional case.
We propose a partition of unity technique to overcome this problem in Banach spaces.
Application is given to a quasiconvex quasioptimization problem through the use of a new existence result for generalized quasivariational inequalities which is based on the Schauder fixed point theorem.
\end{abstract}

\begin{keyword}
Cone upper semicontinuity \sep Normal operator \sep Quasivariational inequality \sep Quasioptimization
\end{keyword}

\end{frontmatter}


\section{Introduction}

The notion of upper semicontinuity seems to be unappropriate for cone-valued maps and hence modified definitions were been introduced and studied \cite{AuHa05,BoCr90,LuPe01}.
The normal cone operator to the adjusted sublevel sets of a quasiconvex function $f$ defined on a Banach space was introduced in \cite{AuHa05} and it was proved to be both quasimonotone and cone upper semicontinuous.
In particular, the authors showed that the normal cone operator admits a locally defined base-valued submap being norm-to-weak$^*$ upper semicontinuous.
In \cite{AuCo13}, when the space is Euclidian, the authors obtained a globally defined upper semicontinuous base-valued submap taking the convex hull of the normalized normal operator which is the intersection of the unit sphere, which is compact, with the normal operator, which is closed.
Since the unit sphere is not compact in the infinite dimensional case, this approach is unsuccessful in a Banach space $X$.

The first aim of this paper is to overcome this problem by using a partition of unity technique.
Theorem \ref{th:main} states the existence of a norm-to-weak$^*$ upper semicontinuous submap $A:X\setminus\argmin f\rightrightarrows X^*$ such that each $A(x)$ is a nonempty weak$^*$ compact convex set not containing the origin which generates the normal cone to the adjusted sublevel set at $x$.
Subsequently, we establish an existence result (Theorem \ref{th:existenceGQVI}) for a generalized quasivariational inequality which improves the famous Tan's result \cite{Ta85}.
Finally, combining both results, we present an application to quasioptimization problems.


In the last part of this introduction, we present some preliminary notions and results.
Let $X$ be a real Banach space with norm $\|\cdot\|$, $X^*$ its topological dual with norm $\|\cdot\|_*$, and $\langle\cdot,\cdot\rangle$ the duality pairing between $X^*$ and $X$.
From now on, unless otherwise indicated, the spaces $X$ and $X^*$ will be equipped by the strong (norm) topology $s$ and the weak$^*$ topology $w^*$, respectively.
The closed unit balls in $X$ and $X^*$ are denoted by $B$ and $B^*$, respectively.
Given a nonempty set $A\subseteq X$ and $x\in X$, $\dist(x,A)=\inf\{\|y-x\|:y\in A\}$ is the distance of $x$ from $A$ and $B(A,r)=\{x\in X:\dist(x,A)\le r\}$ is the neighbourhood of $A$ with radius $r\geq 0$.

A subset $K$ of $X^*$ is a cone if for each $x^*\in K$ and scalar $t>0$, the product $tx^*\in K$ (note that some authors define cone with the scalar $t$ ranging over all non-negative scalars).
Clearly the empty set is a cone.

Let $K\subseteq X^*$ be a cone.
A convex subset $A$ of $K$ is called a base if $K=\{tx^*:t\ge 0,\ x^*\in A\}$ and $0\not\in w^*\mbox{-}\cl A$, where $w^*\mbox{-}\cl$ denotes the closure with respect to the weak$^*$ topology.
Clearly the empty set is a base of the empty cone.
Vice versa, if $K$ admits a nonempty base then $K$ is a convex cone such that $\{0\}\subsetneq K$.
In particular, if the base is compact then $K$ is closed.

The domain and the graph of a set-valued map $\Phi:X\rightrightarrows X^*$ are denoted by $\dom\Phi$ and $\gph\Phi$, respectively.
The map $\Phi$ is norm-to-weak$^*$ upper semicontinuous at $x\in X$ if for every open set $\Omega$ such that $\Phi(x)\subseteq\Omega$, there exists a neighbourhood $U_x$ of $x$ such that $\Phi(x')\subseteq\Omega$, for all $x'\in U_x$.
The map $\Phi$ is norm-to-weak$^*$ closed at $x\in\dom f$ if for each $x^*\in\Phi(x)$ and for each net $\{(x_\alpha,x^*_\alpha)\}$ with $x^*_\alpha\in\Phi(x_\alpha)$ which converges to $(x,x^*)$ in the $s\times w^*$ topology, we have that $x^*\in\Phi(x)$.
The map is norm-to-weak$^*$ closed if its graph is closed with respect to the topology $s\times w^*$.
The Closed Graph Theorem states that a closed-valued map $\Phi$ with values in a compact set is norm-to-weak$^*$ upper semicontinuous if and only if it is norm-to-weak$^*$ closed.

When we are dealing with a cone-valued map $\Phi$, the concept of norm-to-weak$^*$ upper semicontinuity is not appropriate to picture the behaviour of $\Phi$ and it is convenient to slightly alter the definition.
The cone-valued map $\Phi$ is called
\begin{itemize}
\item norm-to-weak$^*$ cone upper semicontinuous at $x\in X$ if for every open cone $\Omega$ such that $\Phi(x)\subseteq \Omega\cup\{0\}$, there exists a neighbourhood $U_x$ of $x$ such that $\Phi(x')\subseteq \Omega\cup\{0\}$, for all $x'\in U_x$;
\item norm-to-weak$^*$ base upper semicontinuous at $x\in X$ if there exist a neighbourhood $U_x$ of $x$ and a set-valued map $A:U_x\rightrightarrows X^*$ such that $A(x')$ is a base of $\Phi(x')$ for each $x'\in U_x$ and $A$ is norm-to-weak$^*$ upper semicontinuous at $x$.
\end{itemize}
Some remarks are needed.
If $\Phi$ is norm-to-weak$^*$ base upper semicontinuous at $x\in X$ then there exists a neighbourhood $U_x$ of $x$ such that $\Phi(x')\ne\{0\}$ for each $x'\in U_x$.
Instead, if $\Phi$ is norm-to-weak$^*$ cone upper semicontinuous at $x\not\in\dom\Phi$ then there exists a neighbourhood $U_x$ of $x$ such that $\Phi(x')\subseteq\{0\}$ for each $x'\in U_x$.
Therefore, if $\Phi$ is norm-to-weak$^*$ cone upper semicontinuous and $\Phi(x)$ admits a base for each $x\in X$ then $\dom\Phi$ is closed.
Moreover the norm-to-weak$^*$ base upper semicontinuity of $\Phi$ at $x$ implies the norm-to-weak$^*$ cone upper semicontinuity at the same point.
The reverse implication holds if $\Phi(x)$ admits a base and $\Phi(x')\ne\{0\}$ for all $x'$ in a suitable neighbourhood of $x$ \cite{AuHa05}.

The norm-to-weak$^*$ cone upper semicontinuity of a map implies its norm-to-weak$^*$ closedness if the map admits a compact base at every point \cite[Proposition 2.3]{BiHaPi22}.
The same proof works for local closedness.
\begin{Theorem}\label{th:localclosedness}
Let $\Phi:X\rightrightarrows X^*$ be a cone-valued map which is norm-to-weak$^*$ cone upper semicontinuous at $x\in\dom\Phi$.
If $\Phi(x)$ has a compact base then $\Phi$ is norm-to-weak$^*$ closed at $x$.
\end{Theorem}

\section{The result}

Let $f:X\rightarrow\R\cup\{+\infty\}$ be an extended-valued function.
Define for any $\lambda\in\R\cup\{+\infty\}$ the sublevel and the strict sublevel set of $f$ at level $\lambda$ by $S_{\lambda}=\{x\in X:f(x)\leq\lambda\}$ and $S_{\lambda}^{<}=\{x\in X:f(x)<\lambda\}$, respectively.
Clearly $S_{\infty}=X$ and $S_{\infty}^{<}=\dom f$.
The function $f$ is quasiconvex if $S_{\lambda}$ is convex for all $\lambda\in\R$.
Now, we recall the notion of adjusted level set introduced in \cite{AuHa05}.
\begin{Definition}
Let $f:X\rightarrow\R\cup\{+\infty\}$ and $x\in X$.
The adjusted sublevel set of $f$ at $x$ is
$$
S^a_f(x)=\left\{
\begin{array}{ll}
S_{f(x)} & \mbox{if } x\in\argmin f\\
S_{f(x)}\cap B(S^<_{f(x)},\rho_x) & \mbox{if } x\notin\argmin f
\end{array}
\right.
$$
where $\rho_x=\dist(x,S^<_{f(x)})$.
\end{Definition}
Note that $S_{f(x)}^{<}\subseteq S^a_f(x)\subseteq S_{f(x)}$ for all $x\in X$; moreover the convexity of the adjusted sublevel sets characterizes the quasiconvexity of the function.
\begin{Theorem}[Proposition 2.4 in \cite{AuHa05}]
The extended-valued function $f$ is quasiconvex if and only if $S^a_f(x)$ is convex, for every $x\in X$.
\end{Theorem}
To any function $f$ we associate the set-valued map $N^a:X\rightrightarrows X^*$ defined by
$$
N^a(x)=\{x^*\in X^*:\langle x^*,y-x\rangle\leq 0,\ \forall y\in S^a_f(x)\}
$$
In \cite[Proposition 3.5]{AuHa05} the authors showed that $N^a$ is norm-to-weak$^*$ base upper semicontinuous under regularity assumptions on $f$.
Combining Theorem \ref{th:localclosedness} and Proposition 3.5 in \cite{AuHa05}, the following result can be easily deduced.
\begin{Corollary}\label{cor:localclosedness}
Let $f$ be quasiconvex and lower semicontinuous at $x\in\dom f\setminus\argmin f$.
If there exists $\lambda<f(x)$ such that $\inte S_{\lambda}\ne\emptyset$ then $N^a$ is closed at $x$.
\end{Corollary}
Such a result has been proved in \cite{AuCo13} and, with weaker assumptions but in a finite dimensional case, in \cite{AlHaSh18}.
Taking advantage of Corollary \ref{cor:localclosedness}, the authors deduce \cite[Proposition 4.4]{AuCo13} the upper semicontinuity of the normalized map $N^a\cap S:\R^n\setminus\argmin f\rightrightarrows B$, being the unit sphere $S$ in $\R^n$ compact.
Moreover, the assumptions in \cite[Proposition 4.4]{AuCo13} guarantee that the convex hull of $N^a\cap S$ is an upper semicontinuous base-valued submap of $N^a$.
Our aim is to extend their result to the infinite dimensional case.
Since the sphere is not weak$^*$ compact in the dual of a Banach space, the previous technique does not work.
\begin{Theorem}\label{th:main}
Let $f:X\rightarrow\R\cup\{+\infty\}$ be proper, quasiconvex and lower semicontinuous.
Assume that for each $x\in X\setminus\argmin f$ there exists $\lambda<f(x)$ such that $\inte S_\lambda\ne\emptyset$.
Then there exists a  norm-to-weak$^*$ upper semicontinuous set-valued map $A:X\setminus\argmin f\rightrightarrows B^*$ such that $A(x)$ is a compact base of $N^a(x)$, for all $x$.
\end{Theorem}
{\bf Proof.}
For the first step of the proof, we argue as in \cite[Lemma 3.6]{AuHa05}.
Let $z\in X\setminus\argmin f$ be fixed.
Choose $z_0\in X$ and $\lambda\in\R$ such that $\lambda<f(z)$ and $z_0\in\inte S^<_\lambda$.
Since $f$ is lower semicontinuous, there exists $\varepsilon>0$ such that
$$z_0+2\varepsilon B\subseteq  S^<_\lambda\subseteq S^<_{f(x)},\qquad \forall x\in z+\varepsilon B$$
Thus, for every $x\in z+\varepsilon B$ and for every
$$
x^*\in N^<(x)=\{x^*\in X^*:\langle x^*,y-x\rangle\leq 0,\ \forall y\in S^<_{f(x)}\}
$$
we obtain the following:
$$\langle x^*,z_0+2\varepsilon u-x\rangle\le 0,\qquad\forall u\in B$$
It follows that
\begin{eqnarray*}
2\varepsilon\|x^*\|_* = 2\varepsilon\sup_{u\in B}\langle x^*,u\rangle & \le & \langle x^*,x-z_0\rangle \\
& = & \langle x^*,z-z_0\rangle+\langle x^*,x-z\rangle \\
& \le & \langle x^*,z-z_0\rangle+\varepsilon \|x^*\|_*
\end{eqnarray*}
Thus,
$$\langle x^*,z-z_0\rangle \ge\varepsilon\|x^*\|_*, \qquad\forall x\in z+\varepsilon B,\ x^*\in N^<(x)$$
Set $H_z=\{x^*\in X^*:\langle x^*,z-z_0\rangle=\varepsilon\}$.
Obviously, for every $x\in z+\varepsilon B$ we have $N^<(x)\cap H_z\subseteq B^*$ and, since $N^a(x)\subseteq N^<(x)$, the set $N^a(x)\cap H_z\subseteq B^*$ is a compact base for the cone $N^a(x)$.
Now, following the proof of \cite[Proposition 3.5]{AuHa05}, we get the  norm-to-weak$^*$ upper semicontinuity of the set-valued map $A_z:z+\varepsilon B\rightrightarrows X^*$ defined by $F_z(x)=N^a(x)\cap H_z$, for all $x\in z+\varepsilon B$.

The last step of the proof consists in finding the selection $A$ as convex combination of the local maps $A_z$ through a partition of unity technique.

Since $X\setminus\argmin f$ is paracompact, there exists a locally finite open covering $\calu=\{U_i:i\in I\}$ where every $U_i\in\calu$ is a subset of some ball $z+\varepsilon B$: let us denote by $A_i$ the map $A_z$ corresponding to the ball $z+\varepsilon B$.
Moreover, there is a partition of unity $\{\lambda_i:i\in I\}$ subordinate to $\calu$ such that each $\lambda_i:X\setminus\argmin f\rightarrow [0,1]$ is continuous, the finite sum $\sum_{i\in I}\lambda_i(y)=1$ for any $y$ and $\lambda_i(y)=0$ for each $y\not\in U_i$.
For every $x\in X\setminus\argmin f$, let $I(x)=\{i\in I:\lambda_i(x)>0\}$, which is nonempty and finite, and define the map
$A:X\setminus\argmin f\rightrightarrows X^*$ as follows
$$
A(x)=\sum_{i\in I(x)}\lambda_i(x)A_i(x)
$$
Clearly $A(x)$ is a compact base of $N^a(x)$, for all $x$.
Moreover, since the values of $A$ are all contained in the compact ball $B^*$, the norm-to-weak$^*$ upper semicontinuity of $A$ is equivalent to prove that the graph of $A$ is closed with respect to the  $s\times w^*$ topology.
Assume that the net $\{x_\alpha\}$ converges to $x$.
Since all the $\lambda_i$ are continuous, it is not restrictive to assume that $I(x)\subseteq I(x_\alpha)$ for all $\alpha$ and we get:
$$
A(x_\alpha) =\sum_{i\in I(x)}\lambda_i(x_\alpha)A_i(x_\alpha)+\sum_{i\in I(x_\alpha)\setminus I(x)}\lambda_i(x_\alpha)A_i(x_\alpha)
$$
Moreover, from the continuity of the functions $\lambda_i$, we deduce
\begin{equation}\label{eq}
\lim_\alpha \sum_{i\in I(x_\alpha)\setminus I(x)}\lambda_i(x_\alpha)=1-\lim_\alpha\sum_{i\in I(x)}\lambda_i(x_\alpha)=0
\end{equation}
Now, let $\{x^*_\alpha\}$ be a net which weakly$^*$ converges to $x^*$ and such that $x^*_\alpha\in A(x_\alpha)$ for any $\alpha$.
Then, there exist $x^*_{i,\alpha}\in A_i(x_\alpha)$ for every $i\in I(x_\alpha)$ such that
\begin{equation}\label{eq:sum}
x^*_\alpha = \sum_{i\in I(x)}\lambda_i(x_\alpha)x^*_{i,\alpha} + \sum_{i\in I(x_\alpha)\setminus I(x)}\lambda_i(x_\alpha)x^*_{i,\alpha}
\end{equation}
The second addend of (\ref{eq:sum}) weakly$^*$ converges to zero since, thanks to (\ref{eq}), it converges to zero in norm
$$
\left\|\sum_{i\in I(x_\alpha)\setminus I(x)}\lambda_i(x_\alpha)x^*_{i,\alpha}\right\|_* \le
\sum_{i\in I(x_\alpha)\setminus I(x)}\lambda_i(x_\alpha)\|x^*_{i,\alpha}\|_* \le
\sum_{i\in I(x_\alpha)\setminus I(x)}\lambda_i(x_\alpha)
$$
On the other hand, without loss of generality, we may assume that $\{x^*_{i,\alpha}\}$ weakly$^*$ converges to some $x^*_i$, for every $i\in I(x)$.
Since $A_i$ has closed graph, we obtain $x^*_i\in A_i(x)$ and $x^*\in A(x)$ follows from (\ref{eq:sum}) taking the weak$^*$ limit.
\hfill$\Box$

\section{An application}

In this section, our aim is to consider a special optimization problem, called quasioptimization problem, and to provide an existence result
for this problem through the study of an associated generalized quasivariational inequality where Theorem \ref{th:main} plays a key role.
We start establishing a new existence result for a generalized quasivariational inequality without requiring any assumption of monotonicity.

Let $C$ be a nonempty subset of $X$ and $T:C\rightrightarrows X^*$ and $K:C\rightrightarrows C$ be two set-valued maps; the generalized quasivariational inequality $GQVI(T,K)$ consists in finding
\begin{equation*}
x\in K(x) \text{ such that } \exists x^*\in T(x)\text{ with }\langle x^*,y-x\rangle\geq 0,\quad\forall y\in K(x)
\end{equation*}
One of the most classic existence results for $GQVI(T,K)$ in the infinite dimensional setting is due to Tan and it was originally stated for locally convex topological vector spaces.
We recall that the set-valued map $K:C\rightrightarrows C$ is said to be lower semicontinuous if for every open set $\Omega$ the lower inverse image $\{x\in C: K(x)\cap\Omega\ne\emptyset\}$ is open in $C$.
Moreover $K$ is called compact if $K(C)$ is contained in a compact subset of $C$.
\begin{Theorem}[Theorem 1 in \cite{Ta85}]\label{th:tan}
Let $C$ be compact and convex and $K$ be closed and lower semicontinuous with nonempty convex values.
Assume that $T$ is norm-to-norm upper semicontinuous with nonempty norm compact convex values, then $GQVI(T,K)$ has a solution.
\end{Theorem}
The existence of solutions for $GQVI(T,K)$ can be obtained with a weaker continuity assumption on $T$ than in Theorem \ref{th:tan} if the space $X$ is normed.
To this purpose, we need to recall the notion of inside point of a convex set that appeared in 1956 in a paper by Michael \cite{Mi56}.
The convex set $S\subseteq C$ is a face of $C$ if $x_1,x_2\in C$, $t\in(0,1)$ and $tx_1+(1-t)x_2\in S$ imply $x_1,x_2\in S$.
Let $\calf_C$ be the (possibly empty) collection of all proper closed faces of $\cl C$, which is the closure of $C$
\begin{Definition}
A point $x\in C$ is an inside point if it is not in any proper closed face of $\cl C$.
Denote by
\begin{displaymath}
I(C)=C\setminus\bigcup_{S\in\calf_C}S
\end{displaymath}
the set of the inside points of $C$.
\end{Definition}
A comparison with other notions of relative interior is given in \cite{CaGi15,CaGi20}.
Thanks to this concept of interior point, we can define the following family of convex sets
\begin{displaymath}
\cald(X)=\{C\subseteq X:C \mbox{ is convex and } I(\cl C)\subseteq C\}
\end{displaymath}
It was proved \cite{Mi56} that $\cald(X)$ contains all the convex sets which are either closed, or with nonempty interior, or finite dimensional.
In particular, when $X$ is finite dimensional the class $\cald(X)$ coincides with the family of all convex sets.
Now we are in position to state and prove our existence result.
Let us denote by $\fix K$ the set of the fixed points of $K$.
\begin{Theorem}\label{th:existenceGQVI}
Let $C$ be convex and $K$ be a compact and lower semicontinuous set-valued map with nonempty values in $\cald(X)$, and $\fix K$ closed.
Assume that $T$ is norm-to-weak$^*$ upper semicontinuous with nonempty weak$^*$ compact convex values, then $GQVI(T,K)$ has a solution.
\end{Theorem}
{\bf Proof.}
Notice that $K$ admits a continuous selection thanks to \cite[Theorem 3.2]{CaGi20}.
Hence the Schauder fixed point theorem as formulated in \cite[Proposition 6.3.2]{GrDu03} guarantees $\fix K\neq\emptyset$.

Let us consider the set-valued map $F:\fix K\rightrightarrows X$ defined as
\begin{displaymath}
F(x)=\bigcap_{x^*\in T(x)}\{y\in X:\langle x^*,y-x\rangle<0\}= \left\{y\in X:\max_{x^*\in T(x)}\langle x^*,y-x\rangle<0\right\}
\end{displaymath}
Clearly, $F$ has convex values.
To prove that $F$ has open graph in $\fix K\times X$, it is sufficient to show that the function $m:\fix K\times X\rightarrow\R$ defined as
\begin{displaymath}
m(x,y)=\max_{x^*\in T(x)}\langle x^*,y-x\rangle
\end{displaymath}
is upper semicontinuous.
First, $\fix K$ is compact since closed subset of the compact set which contains $K(C)$.
From \cite[Lemma 17.8]{AlBo06}, the subset $T(\fix K)$ is weak$^*$ compact; hence, it is norm bounded.
Thanks to \cite[Corollary 6.40]{AlBo06} the duality pairing $\langle\cdot,\cdot\rangle$ restricted to $T(\fix K)\times X$ is jointly continuous, where $X$ has its norm topology and $X^*$ has its weak$^*$ topology; hence, \cite[Lemma 17.30]{AlBo06} guarantees the upper semicontinuity of $m$.

By contradiction, assume that $F(x)\cap K(x)\ne\emptyset$ for all $x\in\fix K$.
Fix $(x_0,y_0)\in\gph K$ and define the map $K_0:C\rightrightarrows C$ as
\begin{displaymath}
K_0(x)=
\left\{\begin{array}{ll}
K(x) & \mbox{ if } x\ne x_0\\
\{y_0\} & \mbox{ if } x=x_0
\end{array}
\right.
\end{displaymath}
$K_0$ is compact and lower semicontinuous, and $K_0(x)\in\cald(X)$ for every $x\in C$.
From \cite[Theorem 3.2]{CaGi20} the map $K_0$ admits a continuous selection, hence $K$ is locally selectionable (see Definition 1.10.1 in \cite{AuCe84}).
From \cite[Proposition 1.10.4]{AuCe84} we deduce that also $F\cap K$ is locally selectionable
and \cite[Proposition 1.10.2]{AuCe84} guarantees that $F\cap K$ has a continuous selection $f:\fix K\rightarrow C$.
Therefore, the set-valued map $\Upsilon:C\rightrightarrows C$ defined as
\begin{displaymath}
    \Upsilon(x)=\left\{\begin{array}{ll}
    K(x) & \mbox{ if } x\notin\fix K\\
    \{f(x)\} & \mbox{ if } x\in\fix K
    \end{array}\right.
\end{displaymath}
is lower semicontinuous \cite[Lemma 2.3]{CaGi20} with values in the class $\cald(X)$.
Hence \cite[Theorem 3.2]{CaGi20} guarantees that $f$ can be extended to a continuous selection $\varphi$ for $\Upsilon$.
The Schauder fixed point theorem guarantees that $\varphi$ has a fixed point, that is, there exists $x\in C$ such that $x=\varphi(x)\in\Upsilon(x)$.
Clearly $x\in\fix K$ and this implies $x=f(x)\in F(x)$ which is absurd.
Therefore, there exists $x\in\fix K$ such that $F(x)\cap K(x)=\emptyset$, that is,
\begin{displaymath}
\min_{y\in K(x)}\max_{x^*\in T(x)}\langle x^*,y-x\rangle\ge 0
\end{displaymath}
Invoking the Sion's minimax theorem we deduce that
\begin{displaymath}
\max_{x^*\in T(x)}\min_{y\in K(x)}\langle x^*,y-x\rangle\ge0
\end{displaymath}
which means that $x$ solves the generalized quasivariational inequality.\hfill$\Box$
\begin{Remark}
Let us compare our result with Theorem \ref{th:tan} due to Tan.
The first difference is about the setting: Tan's result works in a locally convex topological vector space instead Theorem \ref{th:existenceGQVI} is stated in a Banach space.
Nevertheless, the other assumptions of Theorem \ref{th:existenceGQVI} are rather weaker than the ones in Theorem \ref{th:tan}.
Maybe, the most significant improvement consists in requiring the norm-to-weak$^*$ upper semicontinuity of $T$ instead of the stronger norm-to-norm upper semicontinuity.
Moreover, the values of $T$ are assumed weakly$^*$ compact instead of norm compact.
Also the assumptions on $K$ are weaker.
In Theorem \ref{th:tan} the map $K$ is closed, which implies the closedness of $K(x)$, for all $x$.
Conversely, in Theorem \ref{th:existenceGQVI} we require only the closedness of $\fix K$, that is necessary for the closedness of $K$, and $K(x)$ may not be closed but belonging to the class $\cald(X)$ only.
Lastly, we do not assume the compactness of $C$, not even its closedness, but only the fact that $K(C)$ is contained in a compact set.
\end{Remark}
Taking advantage of Theorem \ref{th:existenceGQVI} and the good properties of the normal operator $N_a$, our last aim is to obtain an existence result for a quasioptimization problem through the study of a suitable associated generalized quasivariational inequality.

A quasioptimization problem is an optimization problem in which the constraint set is subject to modifications depending on the considered point.
Given $C\subseteq X$ nonempty, $K:C\rightrightarrows C$ and $f:C\rightarrow\R$, a quasioptimization problem consists in finding
$$
x\in K(x)\text{ such that }f(x)\le f(y),\quad\forall y\in K(x)
$$
Clearly, if $K(x)=C$ for all $x\in C$, quasioptimization problem reduces to a classical optimization problem.

\begin{Theorem}\label{th:existenceQOP}
Let $C$ be convex and $K$ be a compact and lower semicontinuous set-valued map with nonempty values in $\cald(X)$, and $\fix K$ closed.
Assume that $f$ is continuous and quasiconvex, then the quasioptimization problem has a solution.
\end{Theorem}
{\bf Proof.}
Let $T:C\rightrightarrows X^*$ be defined as
$$
T(x)=\left\{
\begin{array}{ll}
B^*  & \mbox{if } x\in\argmin f\\
A(x) & \mbox{if } x\notin\argmin f
\end{array}
\right.
$$
where $A$ is the norm-to-weak$^*$ upper semicontinuous set-valued map obtained in Theorem \ref{th:main}.
Since $\argmin f$ is closed and $A(x)\subseteq B^*$, then $T$ is norm-to-weak$^*$ upper semicontinuous.
In this way, thanks to Theorem \ref{th:existenceGQVI}, it follows that $GQVI(T,K)$ has a solution $x\in C$.
Clearly, if $x\in\argmin f$, then $f(x)\leq f(y)$ for all $y\in K(x)$.
Instead, if $x\notin\argmin f$, then it results that
$$
x^*\in T(x)=A(x)\subseteq N^a(x)\setminus\{0\}
$$
Hence, $x$ is a solution to the generalized variational inequality associated to the operator $N^a\setminus\{0\}$ and the feasible set $K(x)$.
Thanks to \cite[Proposition 3.2]{AuYe06}, the thesis follows.\hfill$\Box$

Theorem \ref{th:existenceQOP} extends Proposition 4.5 in \cite{AuCo13} which is stated in a finite dimensional space and requires also the compactness of $C$ and the closedness of $K$.


\end{document}